\documentclass[11pt]{article}
\usepackage{amsmath,amsthm,amssymb}

\theoremstyle{plain}

 \newtheorem{thm}{Theorem}[section]
 \newtheorem{prop}[thm]{Proposition}
 
 \newtheorem{lem}[thm]{Lemma}

\theoremstyle{definition}

 \newtheorem{dfn}[thm]{Definition}

\theoremstyle{remark}

 \newtheorem{rem}[thm]{Remark}
 \newtheorem{exa}[thm]{Example}
 \newtheorem{prob}[thm]{Problem}
 \newenvironment{pf}{\paragraph{Proof}}{\par\smallskip}
 \numberwithin{equation}{section}

 \DeclareMathOperator{\diag}{diag}
 \DeclareMathOperator{\hcf}{hcf}
 \DeclareMathOperator{\wt}{wt}
 \DeclareMathOperator{\im}{im}
 \DeclareMathOperator{\Gr}{Gr}
 \DeclareMathOperator{\GL}{GL}
 \DeclareMathOperator{\OGr}{OGr}
 
 \DeclareMathOperator{\wGr}{wGr}
 \DeclareMathOperator{\aGr}{aGr}
 \DeclareMathOperator{\aSi}{a\Si}
 \DeclareMathOperator{\wSi}{w\Si}
 \DeclareMathOperator{\aOGr}{aOGr}
 \DeclareMathOperator{\wOGr}{wOGr}
 \DeclareMathOperator{\rank}{rank}
 
 \DeclareMathOperator{\Proj}{Proj}
 \DeclareMathOperator{\SL}{SL}
 \DeclareMathOperator{\SO}{SO}
 \DeclareMathOperator{\Spec}{Spec}
 \DeclareMathOperator{\Spin}{Spin}
 \DeclareMathOperator{\Stab}{Stab}
 \DeclareMathOperator{\Hom}{Hom}
 
 \DeclareMathOperator{\Pf}{Pf}

 \newcommand{\recip}[1]{\frac{1}{#1}}
 \newcommand{\rest}[1]{{}_{|{#1}}}
 \newcommand{\1}{^{-1}}
 \newcommand{\half}{\recip2}
 \newcommand{\tensor}{\otimes}
 \newcommand{\iso}{\cong}
 \newcommand{\into}{\hookrightarrow}
 \newcommand{\onto}{\twoheadrightarrow}
 \newcommand{\ot}{\leftarrow}
 \newcommand{\W}{\bigwedge^2}
 \newcommand{\WV}{\W V}
 \newcommand{\Span}[1]{\left<#1\right>}
 \newcommand{\Oh}{\mathcal O}
 \newcommand{\sC}{\mathcal C}
 \newcommand{\sE}{\mathcal E}
 \newcommand{\sF}{\mathcal F}
 
 \newcommand{\sI}{\mathcal I}
 \newcommand{\sK}{\mathcal K}
 \newcommand{\sL}{\mathcal L}
 \newcommand{\sV}{\mathcal V}
 \newcommand{\Ga}{\Gamma}
 \newcommand{\De}{\Delta}
 \newcommand{\Si}{\Sigma}
 
 \newcommand{\la}{\lambda}
 \newcommand{\al}{\alpha}
 \newcommand{\fie}{\varphi}
 
 \newcommand{\fb}{\mathfrak b}
 \newcommand{\fB}{\mathfrak B}
 \newcommand{\fg}{\mathfrak g}
 \newcommand{\ft}{\mathfrak t}
 \newcommand{\bv}{\mathbf v}
 \newcommand{\bw}{\mathbf w}
 \newcommand{\wF}{\widetilde F}
 \newcommand{\aff}{\mathbb A}
 \newcommand{\PP}{\mathbb P}
 \newcommand{\Z}{\mathbb Z}
 \newcommand{\Q}{\mathbb Q}
 \newcommand{\C}{\mathbb C}
 \newcommand{\Cstar}{\C^\times}
 \newcommand{\myline}{\frac{\ \quad}{\quad}}

\title{Weighted Grassmannians\thanks{For the proceedings of
Paolo Francia memorial conference, Genova, Sep 2001, edited by Mauro
Beltrametti, to appear de Gruyter 2002}}
\author{Alessio Corti \and Miles Reid}
\date{}

\begin{document}
\maketitle

 \begin{center}
 {\em In memory of Paolo Francia}
 \end{center}
 \bigskip

\begin{abstract} Many classes of projective algebraic varieties can be
studied in terms of graded rings. Gorenstein graded rings in small
codimension have been studied recently from an algebraic point of view,
but the geometric meaning of the resulting structures is still relatively
poorly understood. We discuss here the weighted projective analogs of
homogeneous spaces such as the Grassmannian $\Gr(2,5)$ and orthogonal
Grassmannian $\OGr(5,10)$ appearing in Mukai's linear section theorem
for Fano 3-folds, and show how to use these as ambient spaces for weighted
projective constructions. This is a first sketch of a subject that we
expect to have many interesting future applications.
\end{abstract}

\section{Introduction}
We are interested in describing algebraic varieties explicitly in terms of
graded rings and, conversely, in algebraic varieties for which such an
explicit description is possible. Our varieties $X$ always come with a
polarisation $A$, usually the canonical class or an integer submultiple of
it. Our favourites include the following:
\begin{enumerate}
 \renewcommand{\labelenumi}{(\arabic{enumi})}
 \item canonical curves, K3 surfaces, Fano 3-folds. A Fano 3-fold $V$ is
canonically polarised by its anticanonical class $A=-K_V$. We consider K3
surfaces with Du Val singularities polarised by a Weil divisor. A
canonical curve $C$ is a curve of genus $g\ge2$ with its canonical
polarisation by $K_C$; but we will more often be concerned with
subcanonical curves, polarised by a divisor $A$ that is a submultiple
of $K_C=kA$, and variants on orbifold curves also occur naturally (see
Alt{\i}nok, Brown and Reid \cite{ABR}).
\item Regular canonical surfaces, Calabi--Yau 3-folds.
\item Regular canonical 3-folds.
\end{enumerate}
If $X,\Oh(1)$ is a polarised $n$-fold with $K_X=\Oh(k)$, Mukai defines
the {\em coindex} of $X$ to be $n+1+k$; the above are varieties of
coindex~3, 4 and~5.

\begin{rem}
 \begin{enumerate}
 \renewcommand{\labelenumi}{(\alph{enumi})}
 \item Describing a variety explicitly means embedding it into a suitable
ambient space and writing down its equations. This is closely related to
the problem of finding generators and relations for the graded ring
 \[
 R(X,A)=\bigoplus_{n=0}^\infty H^0(X,nA).
 \]
 In all the examples we consider, this is a Gorenstein ring; this property
is one of the most powerful general tools we have in studying $X$ and its
deformations. It seems to us that this point is not adequately appreciated.
 \item Varieties often come in ladders of successive hyperplane sections.
For example, in good situations, a general elephant $S\in|{-}K_V|$ on a
nonsingular Fano \hbox{3-fold} $V$ is a K3 surface polarised by
$A=-K_V\rest S$, and a general $C\in|A|_S$ a canonical curve. Finding the
equations of $V$ is closely related to finding the equations of $S$ or $C$,
and often practically equi\-valent to it.
 \item The natural context to study Fano 3-folds is the Mori category of
projective varieties with terminal singularities. The key examples of these
are cyclic quotient singularities
 \[
 \recip{r}(1,a,-a)\ =\ \C^3/(\Z/r\Z),
 \]
 where the notation signifies that the cyclic group $\Z/r\Z$ acts diagonally
with weights $1,a,-a$, and $\hcf(a,r)=1$. We are thus led to consider K3
surfaces with singularities $\recip{r}(a,-a)$ polarised by ample {\em
Weil\/} divisors, and, one further step down the ladder, {\em orbifold\/}
canonical curves; compare \cite{ABR}.
 \end{enumerate}
\end{rem}

Our original motivation is Mukai's description of a prime Gorenstein Fano
\hbox{3-fold} of genus $7\le g\le10$ as a {\em linear section} of a special
projective homogeneous space, that is, the quotient $G/P$ of a (semisimple)
Lie group $G$ by a maximal parabolic subgroup $P$. For example, consider
$V=\C^{2n}$ endowed with a complex symmetric quadratic form $q$; it is
traditional to take
 \[
 q=\begin{pmatrix} 0 & I \\ I & 0
 \end{pmatrix}, \quad\hbox{where $I=I_{n\times n}$,}
 \]
so that $V=U\oplus U^\vee$ where $U=\Span{e_1,\dots,e_n}$. It is well
known that the space of isotropic $n$-dimensional vector subspaces of $V$
splits into two components for reasons of ``spin''. The O'Grassmann or
orthogonal Grassmann variety $\OGr(n,2n)$ is one connected component; we
take the component containing the reference subspace $U$. It is a
homo\-geneous space for the group $G=\SO(2n,\C)$, and has a natural
Pl\"ucker style spinor embedding in the projective space $\PP(S^+)$ of the
spinor representation $S^+=\bigwedge^{\mathrm{even}}U$. Mukai proves the
following result:

\begin{thm}[Mukai \cite{Mu}] A prime Gorenstein Fano $3$-fold of genus $7$
is a linear section of\/ $\OGr(5,10)$ in its spinor embedding. In other
words, there are $7$ hyperplanes $H_1,\dots,H_7$ of\/ $\PP(S^+)=\PP^{15}$
such that
 \[
 (V_{12}\subset\PP^{8}) \ =\ \OGr(5,10)\cap H_1\cap\cdots\cap H_7.
 \]
\end{thm}

We wanted to see how far these ideas of Mukai generalise. In this note,
we define weighted Grassmann and orthogonal Grassmann varieties, and study
some examples of their linear sections.

 \subsection*{Acknowledgements}
 The treatment of weighted projective homogeneous spaces in
Section~\ref{sec!phom} is based on notes of Ian Grojnowski \cite{G}, whom
we had hoped to involve as coauthor. We have also benefitted from
discussion with Jorge Neves, whose forthcoming paper \cite{N} takes these
ideas further. Gavin Brown and Roberto Pignatelli helped us with computer
algebra calculations.

\section{Weighted $\Gr(2,5)$}\label{sec!Gr}

 \subsection*{The affine Grassmannian $\aGr(2,5)$}
 Weighted versions of a projective homogeneous variety $\Si$ arise on
dividing out the affine cone over $\Si$ by different $\Cstar$ actions.
The construction is particularly transparent for $\Gr(2,5)$. Set $V=\C^5$;
the affine Grassmann variety
 \[
 \aGr(2,5)\subset\WV
 \]
can be defined in any of the following equivalent ways:
 \begin{enumerate}
 \renewcommand{\labelenumi}{(\roman{enumi})}
 \item The variety of skew tensors of rank $\le2$, that is, the image
of $V\times V \to\WV$ given by $(a,b)\mapsto a\wedge b$. In coordinates,
we write
 \begin{equation}
 \W \begin{pmatrix}
 a_1&a_2&a_3&a_4&a_5 \\
 b_1&b_2&b_3&b_4&b_5
 \end{pmatrix}
 =\begin{pmatrix}
 c_{12}&c_{13}&c_{14}&c_{15} \\
 &c_{23}&c_{24}&c_{25} \\
 &&c_{34}&c_{35} \\
 &&&c_{45}
 \end{pmatrix},
 \end{equation}
where $c_{ij}=\det\left|\begin{smallmatrix} a_i&a_j \\ b_i&b_j
\end{smallmatrix} \right|$. Our convention is to write out only the
upper diagonal entries of the $5\times5$ skew matrix $(c_{ij})$.
 \item The closed orbit of the highest weight vector
$e_{12}=(1,0,\dots)\in\WV$ under the action of $\GL(5)=\GL(V)$. In other
words, any tensor of rank~2 is in the $\GL(5)$ orbit of
$(1,0,0,0,0)\wedge (0,1,0,0,0)$.
 \item The quotient of the variety $M(2,5)$ of $2\times5$ matrices by
$\SL(2)$ acting on the left: indeed, the ring of invariant functions is
generated by the $2\times2$ minors $c_{ij}=\det\left|\begin{smallmatrix}
a_i&a_j \\ b_i&b_j \end{smallmatrix} \right|$.
 \item The variety defined by the $4\times4$ Pfaffians of the generic
$5\times5$ skew matrix, that is, the Pl\"ucker equations
 \[
 \Pf_{ijkl}=x_{ij}x_{kl}-x_{ik}x_{jl}+x_{il}x_{jk}=0,
 \]
where $x_{ij}$ for $1\le i<j\le5$ are coordinates on $\WV$. The point is
just that setting the Pfaffians of a skew matrix $(x_{ij})$ equal to zero
enforces rank $\le2$.
 \item\label{it!R} In other words, $\aGr(2,5)$ has affine coordinate ring
 \begin{equation}
 R=\C[\aGr(2,5)]=\C[(x_{ij}]/I,
 \end{equation}
where $I$ is the ideal $I=\Span{\Pf_1,\dots,\Pf_5}$, and $\aGr(2,5)=\Spec R$.
 \end{enumerate}

 \subsection*{Equivariant resolution}
 As a prelude to introducing weights and defining $\wGr$, it is convenient
to explain the symmetry group of $\aGr(2,5)\subset\WV$ and to write its
equations and syzygies in their full symmetry. Under the induced action of
$\GL(V)$ on $\WV$, the scalar matrices $\la\cdot I$ act by $\la^2$.
However, the {\em straight} Grassmannian $\Gr(2,5)\subset\PP^9$ is the
quotient of $\aGr(2,5)$ by $\Cstar$ acting on $\WV$ by overall scalar
multiplication by $\mu\in\Cstar$, and this is not covered by the $\GL(5)$
action; the full symmetry group is thus a double cover of $\GL(V)$ (an
index~2 central extension). Rather than introducing notation for the
double cover, we write $L=\C$ with the usual action of $\Cstar$, and view
the Pl\"ucker embedding as $\aGr(2,5)\into\WV\tensor L$, where $\GL(5)$
acts on the first factor and $\Cstar$ on the second. We also write
 \begin{equation}
 D=\det V\tensor L^2=\bigwedge^5V\tensor L^2,
 \label{eq!D}
 \end{equation}
 a 1-dimensional representation of $\GL(V)\times\Cstar$. It is useful to
bear in mind the straight homogeneous case, when $D$ pushed forward to
$\PP^9$ corresponds to $\Oh(-2)$, and $L$ to $\Oh(-1)$.

 \begin{prop} \label{pr!aGr}
 There are universal maps of vector bundles over the affine space
$\aff^{10}=\WV\tensor L$
 \[
 M\colon \WV\tensor L\to \C \quad \hbox{(that is, $V\tensor L\to V^\vee$)}
 \]
and
 \[
 \Pf=\Pf M\colon V^\vee\tensor D=\bigwedge^4V\tensor L^2\to\C.
 \]
Note that interpreted intrinsically, $\Pf$ is the second wedge of
$M\colon V\tensor L\to V^\vee$.

 Now write $\Oh$ for the structure sheaf of $\aff^{10}$, and
$M\colon\Oh\tensor\WV\tensor L\to \Oh$, etc., for the above universal maps
viewed as sheaf homo\-morphisms. Then the structure sheaf\/ $\Oh_{\aGr}$
of\/ $\aGr(2,5)$ has a $\GL(5)\times\Cstar$ equi\-variant projective
resolution of the form
 \begin{align}
 0\, \ot\, \hbox{} & \Oh \ \xleftarrow{\,\Pf\,} \ \Oh\tensor V^\vee\tensor D
 \ \xleftarrow{\,M\,}\ 
 \Oh\tensor V\tensor L\tensor D
 \ \xleftarrow{\,^t\!\Pf\,}\ \Oh\tensor L\tensor D^2 \ \ot\ 0 \notag \\
 & \,\Big\downarrow \notag\\
 & \!\Oh_{\aGr}. \label{eq!Gr_res}
 \end{align}
 \end{prop}

 \begin{pf} In coordinates $x_{ij}$ on $\WV\tensor L$, the map $M$ is the
generic $5\times5$ skew matrix $(x_{ij})$ and\/ $\Pf=(\Pf_1,\dots,\Pf_5)$
its vector of Pfaffians. Thus (\ref{eq!Gr_res}) follows at once from the
well known fact that the ideal of $\aGr(2,5)$ is generated by the 5
Pfaffians, and $M$ is the matrix of syzygies between them.
 \qed
 \end{pf}

 \begin{rem}\label{rem!eqv}
 Each term in (\ref{eq!Gr_res}) is a $G$-equivariant bundle, where
$G=\GL(5)\times\Cstar$, and the complex gives the projective resolution
of $\Oh_{\aGr}$ in terms of $G$-equivariant vector bundles on the ambient
space $\WV\tensor L$. We write out the definitions for completeness.

Let $G$ be a group and $Y$ a space with a left $G$-action $G\times Y\to
Y$; write $l_g\colon Y\to Y$ for the action of $g\in G$. A $G$-equivariant
sheaf is a sheaf $\sF$ on $Y$, together with isomorphisms
 \begin{equation}
 \al_g\colon\sF\to l_g^*\sF \quad\hbox{satisfying}\quad
 \al_{g_2g_1}=l_{g_2}^*(\al_{g_1})\circ\al_{g_2}.
 \label{eq!cc}
 \end{equation}
 The collection of maps $\{\al_g\}$ is called a {\em $G$-linearisation} or
{\em descent data} for $\sF$; the cocycle condition in (\ref{eq!cc})
ensures that $G$ acts on the pushforward $\pi_*\sF$, where $\pi\colon Y\to
X=Y/G$ is the quotient morphism. A quasicoherent sheaf $\sF$ over an
affine scheme $Y=\Spec A$ is the associated sheaf $\sF=\wF$ for an
$A$-module $F=\Ga(Y,\sF)$; a $G$-equivariant sheaf arises in the same way
from a module $F$ over the twisted group ring $A*G$. That is, $F$ is an
$A$-module with a representation of $G$ such that $g(am)=g(a)g(m)$ for all
$a\in A$ and $m\in F$, where $G$ has the left action on $A$ by
$g(a)=l_{g\1}^\# (a)$.

If $G$ acts freely with quotient $X=Y/G$, taking pushforward and invariant
sections identifies a $G$-equivariant sheaf with a sheaf on $X$; if $G$
has fixed points, the same construction only gives an orbi-sheaf (or a
sheaf on the {\em quotient stack} $[Y/G]$ of which $X$ is the coarse
moduli space). With a little common sense, we can mostly ignore this
point, and pretend that we get a genuine sheaf on the space $Y/G$.
 \end{rem}

 \begin{rem}
 The $\WV\tensor L$ occurring here tells us how to define
$\Gr(2,5)$-bundles over an arbitrary base scheme $S$, more or less as for
conic bundles: choose a rank~5 vector bundle $\sV$, a line bundle $\sL$
and a morphism $\mu\colon\W\sV\tensor\sL\to\Oh_S$, and take the locus
$\rank\mu\le2$ defined by the relative equations $\Pf(\mu)=\W\mu=0$. If
$\sL$ is a square, say $\sL=\sL_0^2$, we can get rid of it by replacing
$\sV\mapsto\sV\otimes\sL_0$.
 \end{rem}

 \subsection*{The definition of $\wGr(2,5)$}
Choosing weights on $\aGr(2,5)$ is equivalent to specifying a 1-parameter
subgroup $\Cstar\into\GL(5)\times\Cstar$. Up to conjugacy, we can
choose it in the maximal torus, that is, diagonal of the form
 \[
 \Bigl(\diag(\la^{w_1},\la^{w_2},\la^{w_3},\la^{w_4},\la^{w_5});\la^u\Bigr)
 \ \subset\ \GL(5)\times\Cstar.
 \]
In order to put weights on $\aGr(2,5)\subset\WV$, we thus specify integer weights
$(w_1,\dots,w_5)$ on $V$, and a separate overall weight $u$ on $\WV$. The
ambient space $\WV$ thus has coordinates
 \[
 x_{ij} \quad\hbox{with}\quad \wt x_{ij}=w_i+w_j+u.
 \]
 Replacing $w_i\mapsto w_i-[\frac u2]$, we can always take $u=0$ or 1. In
fact, for brevity in calculations, we usually use the trick of absorbing
the weight $u$ into the $w_i$ by $w_i\mapsto w_i+\frac u2$, at the cost of
working with half-integers $w_i$. For odd $u$ this is formally incorrect,
but completely harmless, and hardly ever leads to confusion.

\begin{dfn}
 Let $w=(w_1,\dots,w_5)$ and $u$ be weights such that $w_i+w_j+u>0$ for
all $i,j$. We define
 \[
 \wGr(2,5)\ =\ \Bigl(\aGr(2,5)\setminus0\Bigr)/\Cstar,
 \]
where $\Cstar$ acts on $\aGr(2,5)\subset\WV$ by $x_{ij}\mapsto
\la^{w_i+w_j+u}x_{ij}$. Clearly
 \[
 \wGr(2,5)=\Proj R
 \]
where $R=\C[\aGr(2,5)]$ is the affine coordinate ring as in (v) above,
graded by $\wt x_{ij}=w_i+w_j+u$. By definition, $\wGr(2,5)$ comes
with a Pl\"ucker embedding in weighted projective space (wps)
$\PP^9(\{w_i+w_j+u\})$, and is defined by the usual Pl\"ucker equations,
the $4\times4$ Pfaffians of the generic $5\times5$ skew matrix $(x_{ij})$.
 \end{dfn}

 The elementary properties of $\wGr(2,5)$ are easy enough to figure out.
We get affine charts by setting $x_{ij}\ne0$, where (say)
 \[
 x_{12}\ =\ \det \left| \begin{matrix} a_1 & b_1 \\ a_2 & b_2 \end{matrix}
 \right|.
 \]
This chart is the quotient $\C^6/(\Z/\wt x_{12})$ of $\C^6$ by the cyclic
group of order $\wt x_{12}=w_1+w_2+u$ acting on coordinates
$b_3,b_4,b_5,a_3,a_4,a_5$ with weights
 \begin{align*}
 & w_1+w_3+u,\quad w_1+w_4+u,\quad w_1+w_5+u, \\
 & w_2+w_3+u,\quad w_2+w_4+u,\quad w_2+w_5+u.
 \end{align*}
This formula shows the point of our shorthand setting $u=0$, allowing the
$w_i$ to be half-integers. As with weighted projective spaces, we usually
impose ``well formed'' conditions to ensure that the cyclic group acts
effectively and without ramification in codimension~1. We omit the
details, but compare \cite{Fl}, Definition~6.9.

 The Hilbert numerator of $\wGr(2,5)\subset\PP(\{w_i+w_j+u\})$ is
 \[
 \prod_{i,j}(1-t^{w_i+w_j+u})P(t)\ =\ 
 1-\sum_{i=1}^5t^{d-w_i}+\sum_{j=1}^5t^{d+w_j+u}-t^{2d+u},
 \]
 where $d=\sum w_i+2u$. This formula is a numerical version of
(\ref{eq!Gr_res}), and essentially equivalent to it by the splitting
principle. Multiplying by $(1-t)^3$, we deduce that
 \[
 \deg\wGr\ =\ \frac{\sum\binom{d-w_i}3-\sum\binom{d+w_i+u}3+\binom {2d+u}3}
 {\prod (w_i+w_j+u)}\,.
 \]
If $\wGr(2,5)$ is well formed, its canonical class is
$K_{\Gr(2,5)}=\Oh(-2d-u)$. In fact the wps has $K=-\det(\WV\tensor L)$,
which has degree
 \[
 -\sum\wt x_{ij}=-4\sum w_i-10u=-4d-2u,
 \]
and $\wGr(2,5)\subset\PP(\WV\tensor L)$ has the adjunction number
$\deg(L\tensor D^2)=2d+u$ by (\ref{eq!Gr_res}).

 \section*{Tautological sequences}
 Tautological vector bundles over $\aGr(2,5)$ can be discussed in several
ways, parallel to the different treatments of $\aGr(2,5)$. Taking
invariants of the $\Cstar$ action gives rise to tautological
(orbi-)bundles on $\wGr(2,5)$, as in the case of the straight
Grassmannian. These sheaves on $\wGr(2,5)$ can also be understood in terms
of the well known Serre correspondence
 \[
 \sE\mapsto E_*=\bigoplus_{k\ge0}H^0(\sE(k)).
 \]
We describe the Serre module of the tautological bundles explicitly as
modules over the affine coordinate ring of $\aGr$.

 First, $\aGr$ is locally a codimension~3 complete intersection wherever
the matrix of syzygies $M$ has $\rank2$, that is, at every point of
$\aGr\setminus0$. At any such point, we can use two rows of $M$ to express
2 of the 5 Pfaffians as linear combinations of the others, so that the
ideal sheaf $\sI_{\aGr}$ is locally generated by 3 Pfaffians. Thus the
conormal sheaf to $\aGr(2,5)$ is a vector bundle of rank~3 outside the
origin with 5 sections. In more detail, consider (\ref{eq!Gr_res}) as a
resolution of the ideal sheaf $\sI_{\aGr}$:
 \[
 0\ot\sI_{\aGr} \xleftarrow{\Pf} \Oh\tensor V^\vee\tensor D
 \xleftarrow{M} \Oh\tensor V\tensor L\tensor D \ot \cdots.
 \]

Tensoring with $\Oh_{\aGr}=\Oh/\sI_{\aGr}$ gives the exact sequence
 \begin{equation}
 0\ot\sI/\sI^2 \xleftarrow{\Pf} \Oh_{\aGr}\tensor V^\vee\tensor D
 \xleftarrow{M} \Oh_{\aGr}\tensor V\tensor L\tensor D \ot \cdots;
 \label{eq!exa}
 \end{equation}
 twisting back by $D\1$ gives a tautological exact sequence
 \begin{equation}
 0\ot\sF\xleftarrow{\Pf}\Oh_{\aGr}\tensor V^\vee\ot\sE\tensor L\ot0
 \label{eq!F}
 \end{equation}
of vector bundles over $\aGr\setminus0$, where $\sF=\sI/\sI^2\tensor D\1$
and the identification $\ker\Pf=\im M=\sE\tensor L$ is justified below.

Next, $M$ has rank~2 at every point of $\aGr(2,5)\setminus0$, so that if
we restrict the sheaf homomorphism $M\colon\Oh\tensor V\tensor L\tensor
D\to\Oh\tensor V^\vee\tensor D$ to $\aGr(2,5)$, this restriction maps onto
a $\GL(5)\times\Cstar$ equivariant sheaf over $\aGr(2,5)$ that is a rank~2
vector bundle on $\aGr(2,5)\setminus0$. We twist it back by $L\1D\1$ for
convenience, obtaining a second tautological exact sequence:
 \begin{equation}
 0\ot\sE\ot\Oh_{\aGr}\tensor V\ot\sK\ot0.
 \label{eq!K}
 \end{equation}
Here $\sE$ is the same sheaf as in (\ref{eq!F}), up to the indicated
twist, because the sequence in (\ref{eq!exa}) is exact. By playing with
determinant bundles in (\ref{eq!F}) and (\ref{eq!K}) one sees that
$\det\sE=L=\det\sF$ so that $\sE^\vee=\sE\tensor L$, and then the two
sequences are dual, which determines the kernel in (\ref{eq!K}):
 \begin{equation}
 0\ot\sE\ot\Oh_{\aGr}\tensor V\ot \sF^\vee\ot0.
 \label{eq!E}
 \end{equation}

We can concatenate the exact sequences (\ref{eq!E}) and (\ref{eq!F}) to
obtain the following explicit description of the module
$E_*=H^0(\aGr(2,5),\sE)$ over the affine coordinate ring
$R=\C[(x_{ij})]/I=\C[\aGr(2,5)]$: it is generated by 5 sections
$s_1,\dots,s_5$ that one identifies either with the columns of $M$, or
with the 5 columns $s_i=\left(\begin{smallmatrix} a_i\\ b_i
\end{smallmatrix}\right)$ subject to the 10 relations
 \begin{equation}
 x_{ij}s_k - x_{ik}s_j + x_{jk}s_i \quad\hbox{for $1\le i<j<k\le 5$.}
 \label{eq!rels}
 \end{equation}
We can say the same thing in invariant terms by taking global sections in
the exact sequence
 \[
 0\ot\sE\ot\Oh_{\aGr}\tensor V\ot \Oh_{\aGr}\tensor\W V^\vee.
 \]

 Our 4th and final treatment of the bundle $\sE$ is intrinsic and starts
from the model $(\aGr(2,5)\setminus0)=M(2,5)^*/\SL(2)$. Consider the
given representation of $\SL(2)$ on $\C^2$ and the diagonal action of
$\SL(2)$ on the trivial bundle $M(2,5)\times\C^2$; the quotient is the
total space of a rank~2 vector bundle $\sE$ on $\aGr(2,5)\setminus0$. The
sections of $\sE$ are functions $f\colon M(2,5)^*\to\C^2$ that transform
as
 \[
 f(gM)=gf(M) \quad\hbox{for all $g\in\SL(2)$ and $M\in M(2,5)$.}
 \]
 We can identify this bundle $\sE$ with any of the above constructions:
the 5 columns of $M$ give global sections of $\sE$, and they satisfy the
same relations as in (\ref{eq!rels}), leading to the same presentation of
$\sE$ by $\GL(5)\times\Cstar$-equivariant free sheaves on
$\aGr(2,5)\setminus 0$. The advantage of this construction is that, since
it involves the 2-planes parametrised by points of $\aGr(2,5)$, it really
relates to the functor represented by $\aGr(2,5)$, and thus to the
traditional tautological bundle of a Grassmannian.

 \section*{Examples}
 \begin{exa}
 Take $w=(\half,\half,\half,\half,\frac32)$; then $\wGr(2,5)\subset
\PP(1^6,2^4)$ is given by a skew matrix
 \[
 M=\begin{pmatrix} x_{12}&x_{13}&x_{14}&y_1\\
 &x_{23}&x_{24}&y_2\\ &&x_{34}&y_3\\ &&&y_4 \end{pmatrix}
 \quad\hbox{with weights}\quad
 \begin{pmatrix} 1&1&1&2\\ &1&1&2\\ &&1&2\\ &&&2 \end{pmatrix}.
 \]
 This has $5$ Pfaffians of degrees $2,3,3,3,3$. The section
$V=\wGr(2,5)\cap(2)^3$ by 3 general forms of weight~2 is a Fano 3-fold with
 \[
 h^0(V,-K)=6,\quad-K^3=6+\half;
 \]
 that is, $g=4$ and $V$ has in general a singular point $\half(1,1,1)$;
the corresponding family of K3 surfaces is No.~2 in Alt{\i}nok's list,
$\mathtt{Altinok3(2)}$ in the Magma database.

Because the 3 equations of degree~2 are general, they involve 3 of the
weight~2 coordinates $y_1,y_2,y_3$ with nonzero coefficients, and we can
use them to eliminate $y_i$ as generators. Thus we say that $V$ is a {\em
quasilinear section} of $\wGr(2,5)$. By analogy with Mukai's results, we
want to call this a {\em linear section theorem}, but we keep the
``quasi'' for the moment to keep away the unclean spirit.

 The section $S=\wGr(2,5)\cap(2)^4$ by 4 forms of weight~2 has been
studied in detail by Neves \cite{N}; we can use the 4 equations to write
$y_i=q_i(\underline x)$ for $i=1,\dots,4$, giving a canonical surface
$S\subset\PP^5$ with $p_g=6$, $K^2=13$ defined by the Pfaffians of
 \[
 \begin{pmatrix} l_{12} & l_{13} & l_{14} & q_1 \\
 & l_{23} & l_{24} & q_2 \\
 & & l_{34} & q_3 \\
 & & & q_4
 \end{pmatrix}
 \]
with $l_{ij}$ linear and $q_i$ quadratic forms on $\PP(1^6,2^4)$.
Conversely (and slightly more generally), Neves \cite{N} shows that a
surface $S$ with $p_g=6$, $K^2=13$ satisfying appropriate generality
assumptions has a {\em nongeneral} canonical curve $C\in|K_S|$ for which the
restricted linear system splits as $|K_S|\rest C=g^1_6+g^1_7$. Following
Mukai's strategy, Neves shows how to derive the ``tautological'' rank~2
vector bundle $E$ over $S$ and the embedding of $S$ into $\wGr(2,5)$ or a
cone over it from this Brill--Noether data on $C$. The linear entries
$l_{ij}$ of $M$ may be linearly dependent, corresponding to a model of $S$
as a section of a cone over $\wGr(2,5)$.
 \end{exa}

\begin{exa}
 Taking $w=(\half,\half,\half,\frac32,\frac32)$ gives
$\wGr(2,5)\subset\PP(1^3,2^6,3)$ defined by a matrix with weights
 \[
 \begin{pmatrix} 1&1&2&2\\ &1&2&2\\ &&2&2\\ &&&3 \end{pmatrix},
 \]
having Pfaffians of degrees $3,3,4,4,4$.
 \begin{enumerate}
 \renewcommand{\labelenumi}{(\alph{enumi})}
 \item Write $\sC\wGr(2,5)\subset\PP(1^4,2^6,3)$ for the projective cone
over $\wGr(2,5)$; this means that we add one extra variable of degree~1 to
the homogeneous coordinate ring, not involved in any relation. Then a
general quasilinear section $S=\sC\wGr(2,5)\cap(2)^5$ of the cone by $5$
general forms of degree $2$ is a K3 surface with ample Weil divisor $D$
satisfying
 \[
 h^0(S,D)=4,\quad D^2=4+\frac{2}{3}
 \]
 that is, $g=3$ and $S$ has a singular point $\recip3(1,2)$. This family
of K3 surfaces is $\mathtt{Altinok3(3)}$.

The new phenomenon in this example is that $S$ has $h^0(S,D)=4$, so that the
graded ring $R(S,D)$ has 4 generators $x_1,\dots,x_4$ of degree~1. On the
other hand, the matrix only has 3 entries of weight~1, so that not all the
$x_i$ can appear as degree~1 terms. Thus $S$ is obtained from the cone
over $\wGr(2,5)$.

 \item A general quasilinear section $S=\wGr(2,5)\cap (2)^3\cap(3)$ is a
K3 with ample $D$ such that
 \[
 h^0(S,D)=3 \quad\hbox{and}\quad D^2=2+3\times\half
 \]
 that is, $g=2$ and $S$ has 3 singular points $\half(1,1)$; this is
$\mathtt{Altinok3(5)}$.
 \end{enumerate}
\end{exa}
The following is due to Selma Alt\i nok.
 \begin{thm}[Alt{\i}nok \cite{Al}] There are precisely $69$ families of
K3 surfaces with cyclic singularities $\recip{r}(a,-a)$ whose general
element is a codimension~$3$ subvariety in weighted projective space given
by the $4\times4$ Pfaffians of a skew $5\times5$ matrix.
 \end{thm}
The next result is a nice structural description of these surfaces;
unfortunately, we don't know how to prove it in an entirely conceptual way.
\begin{prop}
 All K3 surfaces of Alt\i nok are quasilinear sections of a weighted
Grassmannian $\wGr(2,5)$ or a cone over $\wGr(2,5)$.
\end{prop}

\begin{pf} Ultimately, this is based on a case by case check against
Alt\i nok's list. By Buchsbaum--Eisenbud, $S\subset \PP(a_1,\dots,a_6)$
is defined (scheme theoretically) by the $4\times4$ Pfaffians of a
$5\times5$ skew matrix
 \[
 \begin{pmatrix} f_{12} & f_{13} & f_{14} & f_{15} \\
 & f_{23} & f_{24} & f_{25} \\
 & & f_{34} & f_{35} \\
 & & & f_{45}
 \end{pmatrix}
 \]
 where the entry $f_{ij}(y_1, \dots, y_6)$ is a weighted homogeneous form
of degree $d_{ij}$ in the coordinates $y_i$ (the condition for $S$ to be a
K3 implies that, if $b_i=\deg\Pf_i$ is the degree of the $i$th Pfaffian,
then $\sum b_i=2\sum a_i$). Now an easy combinatorial argument shows that
the 5 Pfaffians are weighted homogeneous, if and only if $d_{ij}=w_i+w_j$
for some $w_i$, $i=1,\dots,5$. The idea is to map $S$ to the weighted
Grassmannian $\wGr(2,5)$ with weights $w=(w_1,\dots,w_5)$, immersed in
$\PP(x_{ij})$, by setting
 \[
 f_{ij}=x_{ij}
 \]
 and check that this is an embedding and maps to a quasilinear section. As
far as we can tell, these must be checked explicitly on each of the 69
families of Alt\i nok. Intuitively, the key point is that, for $S$ to be
a K3, the degrees $f_{ij}$ must be ``small''. Slightly more precisely, the
formula $\sum b_i=2\sum a_i$ implies in practice that many of the $d_{ij}$
equal an $a_i$, which is to say that $f_{ij}$ is linear in one of the
variables.
 \qed \end{pf}

 \begin{rem}
If $S\subset\wGr(2,5)$ is a K3 quasilinear section, it is tempting to try
to reconstruct the embedding from intrinsic data on $S$, by analogy with
Mukai's constructions and Neves \cite{N}. It is easy to see that $E_{|S}$
is a rigid simple vector bundle, hence stable and uniquely characterised
by its Chern classes and local nature at the singularities. We know many
ad hoc constructions but no unified way to produce the bundle directly on
$S$, and no a priori reason why it must exist. In the case of an
(orbifold) canonical curve $C$, the vector bundles $E_{|C}$ arising from
embeddings in $\wGr(2,5)$ are often interesting and rather exceptional
from the point of view of higher rank Brill--Noether theory.
\end{rem}

\section{Weighted homogeneous spaces}\label{sec!phom}

This section is based on a close reading of part of Ian Grojnowski's
notes \cite{G}. We give the general definition of weighted projective
homogeneous spaces under an algebraic group $G$ and describe an explicit
atlas of coordinate charts on them. A homogeneous variety that is
projective is of course homo\-geneous under a semi\-simple group $G$;
however, weighted homogeneous spaces always involve central extensions, as
we saw with $\SL(5)$, $\GL(5)$ and $\GL(5)\times\Cstar$ in the preceding
section. Thus we work from now on with a reductive group $G$.

\subsection*{Notation}
Let $G$ be a reductive complex algebraic group. We fix a maximal torus and
Borel subgroup $T\subset B\subset G$ and write $\fB=G/B$ for the maximal
flag variety. Let $X=\Hom(T,\Cstar)$ be the lattice of weights (or
characters), and $Y=\Hom(\Cstar,T)$ the dual lattice of 1-parameter
subgroups, with the perfect pairing $\Span{\ ,\,}\colon X\times Y\to\Z$.

Recall that the {\em roots} of $G$ are defined as the weights of $T$
appearing in the adjoint action of $T$ on the Lie algebra $\fg$. We write
$\De\subset X$ for the set of these. A root $\al\in\De$ determines an
involution of the maximal torus $T$, and hence a reflection $r_\al$ of
$X$; these reflections generate the Weyl group $W(G)$. The {\em negative
roots} $-\De_+$ are the roots appearing in $\fb/\ft$. Let $S\subset\De_+$
be the set of simple roots.

\subsection*{Projective homogeneous spaces and parabolic subgroups} A
projective homo\-geneous space under $G$ is a quotient space $\Si=G/P$ by
a parabolic subgroup $P$. Every such is conjugate to a {\em standard}
parabolic subgroup, that is, one containing $B$. A standard parabolic
subgroup $P$ corresponds to the subset of simple roots
 \[
 I=\bigl\{\al\in S\bigm|r_\al(B)\subset P\}\subset S.
 \]
We recover $P$ as follows: let $W_I\subset W(G)$ be the subgroup generated
by $r_\al$ for $\al\in I$; then $P=P_I=BW_IB$. We write $\Si_I=G/P_I$ for
the corresponding projective homogeneous space (or generalised flag
variety).

\subsection*{Dominant weights}
A weight $\chi\in X$ extends to a unique character $B\to\C^*$, and hence
gives rise to a line bundle $\Oh(\chi)$ on $\fB$. A weight $\chi$ is {\em
dominant} if $V_\chi=H^0(\fB,\Oh(\chi))\ne0$; thus the cone $X^+$ of
dominant weights is the effective cone of $\fB$. If $\chi$ is a dominant
weight, $V_\chi$ is an irreducible representation of $G$ with highest
weight vector $v_\chi$. The linear system $|V_\chi|$ is free, and defines
an equivariant morphism $G\to\PP(V_\chi)$ whose image $\Si=G\cdot \C
v_\chi$ is the orbit of the highest weight line $\C v_\chi$. Therefore
$\Si=G/P$ is a projective homogeneous space, where $P=\Stab(\C v_\chi)$ is
a parabolic subgroup. Then $P=P_I$ as above and
$\Si=\Si_I\subset\PP(V_\chi)$ is a {\em generalised Pl\"ucker embedding}.

\subsection*{Definition of weighted homogeneous spaces}
Let $\rho\in Y=\Hom(\Cstar,G)$ be a 1-parameter subgroup, and $u\ge0$ an
overall weight. We use $\rho$ and $u$ to make $V_\chi$ into a
representation of $G\times\Cstar$, with the second factor acting by
 \[
 \la\colon v\mapsto\la^u\rho(\la)\cdot v.
 \]
We assume from now on that this action has only positive weights.

\begin{rem} This never happens if, say, $G$ is semisimple and $u=0$.
However, we can always make it happen by taking $u$ large enough; more
precisely, the weights are all positive if and only if
 \[
 N_w=\Span{\chi,w\rho}+u>0 \quad \hbox{for every $w\in W(G)$.}
 \]
Here we assume that this condition is satisfied.
 \end{rem}

 Then the quotient $\PP(V_\chi)(\rho,u)=(V_\chi\setminus0)/\Cstar$ is a
weighted projective space. {From} the description $\aSi_I\subset
V_\chi=G\cdot\C v_\chi$, we see that $\aSi_I$ is invariant under the
$\Cstar$-action.

\begin{dfn} The weighted homogeneous variety associated to this data is
the quotient
 \[
 \wSi_I\ =\ \Bigl(\aSi_I\setminus0\Bigr)/\Cstar\subset\PP(V_\chi)(\rho,u).
 \]
To stress the choices of the data $\chi,\rho,u$, we write
$\wSi_I=\wSi_I(\chi,\rho,u)$.
\end{dfn}

 \begin{lem}
 We have $\wSi_I=\wSi_I(\chi,\rho,u)=\wSi_I(\chi,w\rho,u)$ for all $w\in
W(G)$.
 \end{lem}

\begin{pf} Almost obvious, but see the explicit coordinatisation given
below. \qed
\end{pf}

\subsection*{Coordinate charts}
We write down explicit $T$-invariant coordinate charts on weighted
homogeneous varieties as quotients of affine spaces by a cyclic group.
This explicit coordinate atlas is useful in studying various properties
of $\wSi_I$.

Let $U^-$ be the unipotent radical of the opposite Borel subgroup $B^-$.
Choose a $T$-equivariant isomorphism $\C^{\De_+}\iso U^-$, where $T$ acts
on $\C^{\De_+}$ by $x\cdot s_\al=\al(x^{-1})s_\al$. Thus a 1-parameter
subgroup $\rho\colon\Cstar\to G$ gives rise to an action of $\Cstar$ on
$\C^{\De_+}$ by
 \[
 \la\cdot s_\al=\la^{-\Span{\al,\rho}}s_\al.
 \]

As a warm up, we start with the maximal flag variety $\fB=G/B$. Then for
each $w\in W(G)$, the image of $wU^-v_\chi$ in $\PP(V_\chi)$ is an open set
of $\fB$, isomorphic to the affine space $\C^{\De_+}$ with $T$-action
twisted by $w^{-1}$. Moreover, the union of these open sets is all of
$\fB$. Thus we get a covering of the weighted flag variety $\wSi$ by
$|W(G)|$ open subsets, each isomorphic to $\C^{\De_+}/\mu_{N_w}$, where
$N_w=u+\Span{\chi,w\mu}$ as before, and $\la\in\mu_{N_w}$ acts by
 \[
 \la\cdot s_\al=\la^{-\Span{\al,w\rho}}.
 \]

We now treat the general case of $\Si_I= G/P_I$, for $I\subset S$. Write
$\De^I$ for the roots that can be written as linear combinations of the
roots in $I$, and $\De^I_+=\De^I\cap\De_+$. We set
 \[
 U^-=U^J\times U_J^- \quad\hbox{where}\quad U_J^-\iso\C^{\De_-^I}
 \quad\hbox{and}\quad U^J\iso\C^{\De_+\setminus \De_+^I}
 \]
($T$-equivariantly). Then the weighted homogeneous space
$\wSi_I(\chi,\rho,k)$ admits a cover by $|W(G)/W_I|$ open charts, each a
cyclic quotient of affine space. The chart corresponding to $w$
is the image of $wU^J(v_\chi)$ in $\PP(V_\chi)(\rho,u)$; it is isomorphic
to $U^J/\mu_{N_w}$ where $\la\in\rho_{N_w}$ acts by
 \[
 \la\cdot s_\al=\la^{-\Span{\al,u\rho}}.
 \]

 \begin{prob}
 As with weighted projective spaces, to use weighted homo\-geneous spaces
$\wSi$ as ambient spaces in which to construct varieties, we need to study
questions such as when a subvariety $X\subset\Si$ is well formed (that is,
no orbifold behaviour in codimension~0 or~1, no quasireflections), or
quasismooth (that is, the affine cone over $X$ is nonsingular); for $\wSi$
itself, it seems reasonable to expect that the $\Cstar$ action on
$\PP(V_\chi)$ is well formed if and only if its action on $\aSi_I$ is. By
analogy with the toric case, there must be straightforward adjunction
formulas for the canonical class of weighted homogeneous spaces $\wSi$,
together with criteria to determine whether the affine cone $\aSi$ is
Gorenstein or $\Q$-Gorenstein, and has terminal or canonical
singularities. The results of the preceding section on $\Gr(2,5)$ raise
the interesting question of writing down the projective resolution of
$\aSi_I\subset V_\chi$ in equivariant terms; Lascoux \cite{La} has related
results in some important cases that might serve as a model. Since $\Si_I$
has the status of a generalised flag variety, it is also interesting to
study the corresponding tautological structures over $\aSi_I$.
 \end{prob}

 \section{Weighted orthogonal Grassmannian $\OGr(5,10)$}
As in the introduction, let $V=\C^{10}$ with a nondegenerate quadratic
form $q$; a change of basis puts $q$ in the normal form
 \[
 q=\begin{pmatrix} 0 & I \\ I & 0
 \end{pmatrix}, \quad
 \hbox{that is,} \quad \hbox{$V=U\oplus U^\vee$, where
 $U=\Span{e_1,\dots,e_5}$.}
 \]
 We write $f_1,\dots,f_5$ for the dual basis of $U^\vee$. A vector
subspace $F\subset V$ is {\em isotropic} if $q$ is identically zero on
$F$. For example, $U$ is an isotropic 5-space. Since $q$ is nondegenerate,
it is clear that an isotropic subspace $F\subset V$ has dimension $\le5$.
We say that a maximal isotropic subspace is a {\em generator} of $q$, or
of the quadric hypersurface $Q:(q=0)\subset\PP(V)$. The parity $\dim
F_\la\cap U$ mod~2 is known to be locally constant in a continuous family
of generators $F_\la$. Thus parity splits the generators into two
connected components. We choose the component containing the reference
subspace $U$. Thus we define the {\em orthogonal Grassmann} or {\em
O'Grassmann variety} $\OGr(5,10)$ by
 \[
 \OGr(5,10)=\left\{ F \in\Gr(5,V) \left|
 \begin{array}{l}
 \hbox{$F$ is isotropic for $q$} \\
 \hbox{and $\dim F\cap U$ is odd}
 \end{array} \right.\right\}
 \]

\subsection*{The Weyl group $W(D_5)$} The study of the algebraic group
$\SO(10,\C)$, its double cover $\Spin(10)\to\SO(10,\C)$, and their
representations is governed by the Weyl group $W(D_5)$, which acts as a
permutation group on every combinatoric set in the theory. We are
particularly interested in two permutation representations of $W(D_5)$
that base respectively the given representation $V$ of $\SO(10)$, and the
space of spinor $S^+=\bigwedge^{\mathrm{even}}U$, which is a
representation of $\Spin(10)$.

The given representation $V=U\oplus U^\vee$ has basis
$e_1,\dots,e_5,f_1,\dots,f_5$, and $W(D_5)$ acts by permuting the indices
$\{1,\dots,5\}$ on the $e_i$ and $f_i$ simultaneously, and by swapping
evenly many $e_i$ with $f_i$. For example, the permutations
 \[
 (e_1f_1)(e_2f_2) \quad\hbox{and}\quad (e_1f_1)(e_2f_2)(e_3f_3)(e_4f_4)
 \]
 are elements of $W(D_5)$. One checks that in this permutation
representation, $W(D_5)$ is the Coxeter group generated by the 5
involutions
 \[
 \renewcommand{\arraycolsep}{0.1cm}
 \renewcommand{\arraystretch}{1.4}
 \begin{matrix}
 (12) &\myline& (23) &\myline& (34) &\myline& (45) \\
 &&&& \big| \\
 &&&& \kern-1cm (e_4f_5)(e_5f_4) \kern-1cm
 \end{matrix}
 \]
with the Coxeter relations indicated by the Dynkin diagram.

The spinor representation $S^+$ has basis the 16 nodes of the graph
 \[
 \Ga\ =\ \hbox{5-cube modulo antipodal identification.}
 \]
The action of $W(D_5)$ on $\Ga$ has 5 involutions parallel to the facets
of the 5-cube, whose product is the antipodal involution, and thus acts
trivially on $\Ga$. These define a normal subgroup
$(\Z/2)^5/(\diag)\triangleleft W(D_5)$, the quotient by which is the
symmetric group $S_5$ permuting the 5 orthogonal directions of the 5-cube.
Thus $W(D_5)$ is the extension $(\Z/2)^4\triangleleft W(D_5)\onto S_5$.

To introduce notation for the nodes of $\Ga$, we break the symmetry by
choosing a preferred node $x=x_\emptyset\in\Ga$ and an order $1,2,3,4,5$
on the 5 edges $xx_1,\dots,xx_5$ out of $x$. Then $\Ga$ consists of $x_I$,
where $I\subset\{1,2,3,4,5\}$, and $x_I=x_{\sC I}$ (where $\sC I$ is the
set complement $\sC I=\{1,2,3,4,5\}\setminus I)$. The {\em short}
representatives are $x,x_i,x_{ij}$, with $x=x_{\emptyset}=x_{12345}$,
$x_1=x_{2345}$, $x_{12}=x_{345}$, etc. We use this below to work out the
equations and syzygies of $\OGr(5,10)$.

The symmetry here is the same as that of the 16 lines on the del Pezzo
surface of degree~4, see Reid \cite{R}.

 \subsection*{Notation}
 We take the construction $S^+=\C\oplus \W U\oplus\bigwedge^4U$ as the
definition of $S^+$, without attempting to deal with it intrinsically
(which can be done in terms of the even Clifford algebra). This
construction depends on the choice of $U$ or of the decomposition
$V=U\oplus U^\vee$, and $S^+$ is a representation of the double cover
$\Spin(10)$, not of $\SO(10)$ itself.

We write $(e,M,P)\in S^+$ for an element of $S^+$, where $e\in\C$,
$M=(x_{ij})$ is a skew $5\times5$ matrix and $P$ a $5\times1$ column
vector. If $M=(x_{ij})$ is a skew $5\times5$ matrix then $\Pf M$ is the
column vector of its Pfaffians, that is,
 \[
 \Pf M=
 \begin{pmatrix}
 x_{23}x_{45}-x_{24}x_{35}+x_{25}x_{34} \\
 -x_{13}x_{45}+x_{14}x_{35}-x_{15}x_{34} \\
 x_{12}x_{45}-x_{14}x_{25}+x_{15}x_{24} \\
 -x_{12}x_{35}+x_{13}x_{25}-x_{15}x_{23} \\
 x_{12}x_{34}-x_{13}x_{24}+x_{14}x_{23} \\
 \end{pmatrix}.
 \]

 \subsection*{Affine cover of $\OGr(5,10)$}
 Since $V=U\oplus U^\vee$, a 5-plane $F\subset V$ near $U$ is the graph of
a linear map $\fie\colon U\to U^\vee$, so that $U\in\Gr(5,10)$ has an
affine neighbourhood parametrised by $\Hom(U,U^\vee)$: in other words, $F$
has a basis of 5 vectors in $V=\C^{10}$ that we can write as a matrix
$(I,M)$ with $I=I_5$ and $M$ a $5\times5$ matrix. One sees that $F$ is
isotropic for $q$ if and only if the linear map $\fie$ or the matrix $M$
is skew. Thus an affine neighbourhood of $U$ in $\OGr(5,10)$ is given by
$(I,M)$ with $M$ a skew $5\times5$ matrix.

There are 16 standard affine pieces of $\OGr(5,10)$. Each is obtained from
this one by acting on the basis of $V$ by a permutation of $W(D_5)$. That
is, take matrices $(I,M)$ with $M$ skew, and swap evenly many of the first
5 columns with the corresponding columns from the last 5.

 \subsection*{The spinor embedding $\aOGr(5,10)\subset S^+$}
 Corresponding to the different treatment of tensors of rank~2 in $\WV$ in
Section~\ref{sec!Gr}, we can write down 4 characterisations of {\em
simple} spinors. For a spinor $s\in S^+$, we have the following equivalent
conditions:

 \begin{enumerate}
 \renewcommand{\labelenumi}{(\roman{enumi})}
 \item Explicit: $s$ is in the $W(D_5)$-orbit of a spinor of the form
 \[
 e(1,M,\Pf M) \quad\hbox{with $e\in\C$ and $M$ a skew $5\times5$ matrix}.
 \]

 \item Orbit of highest weight vector: $s$ is in the $\Spin(10)$-orbit of
the spinor $(1,0,\dots,0)\in S^+$.

 \item Quotient by $\SL(5)$: consider all $5\times10$ matrices $N$ whose
rows base a generator $\Pi\in\OGr(5,10)$ of $q$, and take the quotient by
$\SL(5)$ acting by left multi\-plication. To explain this briefly, suppose
that $N=(A,B)$ with nonsingular first $5\times5$ block $A$. Then up to the
$\SL(5)$ action, $N$ is of the form $e(I,M)$. The affine embedding into
$S^+$ is then given by $e(1,M,\Pf M)$. Every other $N$ is of this form up
to the action of the Weyl group $W(D_5)$.

 \item Equations: $s=(e,M,P)\in S^+$ satisfies the 10 equations
 \[
 eP=\Pf M,\quad MP=0
 \]
 (see below). The first set of 5 equations with $e\ne0$ describes the
embedding of the first affine open in spinor space. On the other hand, as
we discuss next, this set of equations is $W(D_5)$-invariant.
 \end{enumerate}

 \begin{rem}
 We use the following point of view in (iii): it is well known that the
spinor embedding $\OGr\into\PP(S^+)$ is the Veronese square root of the
Pl\"ucker embedding $\OGr\into\Gr(5,10)\into\PP(\bigwedge^5\C^{10})$. In
other words, up to a straightforward (!) change of coordinates, the set of
$5\times5$ minors of $N$ is the second symmetric power of the set of
spinor coordinate functions $e,x_{ij},\Pf_k$.
 \end{rem}

 \subsection*{Equations of $\aOGr(5,10)$}

As described above, $S^+$ has a basis indexed by the graph $\Ga$. A pair
$x_I,x_J$ is an {\em edge} of $\Ga$ (that is, $x_I$ is joined to $x_J$) if
and only if $I$ and $J$ or $I$ and $\sC J$ differ by one element. Because
of this definition, edges of $\Ga$ fall into 5 sets of 8 parallel edges,
with {\em directions} given by adding the same $i$: for example, the 8
edges
 \[
 xx_1, \quad x_ix_{1i}, \quad x_{ij}x_{kl} \hbox{ with
$\{i,j,k,l\}=\{2,3,4,5\}$}
 \]
are all of the form $x_Ix_{I+1}$, so in the 1 direction. Two edges of $\Ga$
are {\em remote} if no edge of $\Ga$ joins either end of one to either end
of the other. Any two parallel edges either form two sides of a square, or
are remote. Each set of 8 parallel edges breaks up into two {\em remote
quads}. For example, the 8 edges in the 1 direction give
 \[
 xx_1,x_{23}x_{45},x_{24}x_{35},x_{25}x_{34}
 \quad\hbox{and}\quad
 x_1x_{12},x_1x_{13},x_1x_{14},x_1x_{15},
 \]
The 10 equations of $\OGr(5,10)$ in (iv) are sums of these quads with
appropriate choice of signs:
 \[
 xx_1-x_{23}x_{45}+x_{24}x_{35}-x_{25}x_{34}=0
 \quad\hbox{and}\quad
 x_1x_{12}-x_1x_{13}+x_1x_{14}-x_1x_{15}=0,
 \]
and permutations.

The 10 equations {\em centred at $x$} are written in terms of the matrices
 \[
M= \begin{pmatrix}
 x_{12} & x_{13} & x_{14} & x_{15} \\
 & x_{23} & x_{24} & x_{25} \\
 && x_{34} & x_{35} \\
 \multicolumn{2}{c}{-\hbox{sym}}&& x_{45}
 \end{pmatrix} \quad\hbox{and}\quad
 \bv =\begin{pmatrix} x_1 \\ x_2 \\ x_3 \\ x_4 \\ x_5
 \end{pmatrix}.
 \]
They take the form
 \[
 \begin{pmatrix}
 N_1 \\ N_2 \\ N_3 \\ N_4 \\N_5
 \end{pmatrix}
 =
 x\bv-\Pf M = 0 \quad\hbox{and}\quad
 \begin{pmatrix}
 N_{-1} \\ N_{-2} \\ N_{-3} \\ N_{-4} \\N_{-5}
 \end{pmatrix}
 = M\bv = 0.
 \]

The 16 first syzygies are also indexed by the 16 vertices of $\Ga$. Each
is a 5 term syzygy involving the 5 neighbouring vertices:
{\small
 \[
 \renewcommand{\arraycolsep}{0.15em}
 \kern-.8cm
 \begin{array}{cccccccccccccccc}
x& x_1& x_2& x_3& x_4& x_5& x_{12}& x_{13}& x_{14}& x_{15}&
x_{23}& x_{24}& x_{25}& x_{34}& x_{34}& x_{45} \\[6pt]
\hline\\ [-8pt]
0& 0& x_{12}& x_{13}& x_{14}& x_{15}& x_2& x_3& x_4& x_5& 0&0&0&0&0&0 \\
0& -x_{12}& 0& x_{23}& x_{24}& x_{25}&-x_1& 0& 0& 0& x_3&x_4&x_5&0&0&0\\
0& -x_{13}&-x_{23}& 0& x_{34}& x_{35}& 0& -x_1& 0& 0&-x_2&0&0&x_4&x_5&0\\
0& -x_{14}&-x_{24}&-x_{34}& 0& x_{45}& 0& 0& -x_1& 0&0&-x_2&0&-x_3&0&x_5\\
0& -x_{15}&-x_{25}&-x_{35}&-x_{45}& 0& 0& 0& 0& -x_1&0&0&-x_2&0&-x_3&-x_4\\
x_1& x& 0& 0& 0& 0& 0& 0& 0& 0&-x_{45}&x_{35}&-x_{34}&-x_{25}&x_{24}&-x_{23}\\
x_2& 0& x& 0& 0& 0& 0& x_{45}&-x_{35}& x_{34}&0&0&0&x_{15}&-x_{14}&x_{13}\\
x_3& 0& 0& x& 0& 0& -x_{45}& 0& x_{25}&-x_{24}&0&-x_{15}&x_{14}&0&0&-x_{12}\\
x_4& 0& 0& 0& x& 0& x_{35}& -x_{25}& 0& x_{23}&x_{15}&0&-x_{13}&0&x_{12}&0\\
x_5& 0& 0& 0& 0& x& -x_{34}& x_{24}&-x_{23}& 0&-x_{14}&x_{13}&0&-x_{12}&0&0
 \end{array}
 \kern-.8cm
 \]
}

The 2nd syzygies are likewise indexed by the 16 monomials; they form a
$16\times16$ symmetric matrix with typical columns
 \[
 \kern-.5cm
 \begin{array}{l}
 S(x) \\
x^2\\
xx_1-2N_1 \\
xx_2-2N_2 \\
xx_3-2N_3 \\
xx_4-2N_4 \\
xx_5-2N_5 \\
xx_{12}\\
xx_{13}\\
xx_{14}\\
xx_{15}\\
xx_{23}\\
xx_{24}\\
xx_{25}\\
xx_{34}\\
xx_{35}\\
xx_{45}
 \end{array} \qquad
 \begin{array}{l}
 S(x_1) \\
xx_1-2N_1\\
x_1^2\\
x_1x_2\\
x_1x_3\\
x_1x_4\\
x_1x_5\\
x_1x_{12}+2N_{-2}\\
x_1x_{13}+2N_{-3}\\
x_1x_{14}+2N_{-4}\\
x_1x_{15}+2N_{-5}\\
x_1x_{23}\\
x_1x_{24}\\
x_1x_{25}\\
x_1x_{34}\\
x_1x_{35}\\
x_1x_{45}
 \end{array} \qquad
 \begin{array}{l}
 S(x_{12}) \\
xx_{12}\\
x_1x_{12}+2N_{-2}\\
x_2x_{12}-2N_{-1}\\
x_3x_{12}\\
x_4x_{12}\\
x_5x_{12}\\
x_{12}^2\\
x_{12}x_{13}\\
x_{12}x_{14}\\
x_{12}x_{15}\\
x_{12}x_{23}\\
x_{12}x_{24}\\
x_{12}x_{25}\\
x_{12}x_{34}+2N_5\\
x_{12}x_{35}-2N_4\\
x_{12}x_{45}+2N_3\\
 \end{array}
 \]

\subsection*{Numerology} {From} this we get the following numerology and
representation theory. Write $U=\C^5$ with weights $w_1,\dots,w_5$, where
$\SO(10)$ acts on $V=U\oplus U^\vee$. As in Section~\ref{sec!phom}, to
ensure that all the weights are positive, we introduce a further overall
weight $u$ on $S^+$. To keep track of this, we introduce the bigger group
$G=\Spin(10)\times\Cstar$, and replace $S^+$ by $S^+\tensor L$, where
$\Spin(10)$ acts on the first factor in the usual way, and $\Cstar$ acts
on $L=\C$ with weight $u$. Set
 \[
 S^+=\bigwedge^{\mathrm{even}}U\tensor L=
 (\C\oplus\bigwedge^2 U\oplus\bigwedge^4 U) \tensor L
 \]
for the 16 dimensional spinor space. The only representations we need are
$V,S^+$, its dual $S^-$ and their twists by line bundles. By analogy with
Proposition~\ref{pr!aGr} and (\ref{eq!D}), we define\footnote{As in
(\ref{eq!D}--\ref{eq!Gr_res}), this move cleans up the formulas below in a
most miraculously way. However, to be quite honest, at the time of writing,
we have absolutely no idea what representation $D$ is, or for what group it
is supposed to be an equivariant line bundle on $\C^{16}=S^+$. Cf.\
Problem~\ref{prob!OGr}.}
 \[
 D=\bigwedge^5U\tensor L^2, \quad s=\sum w_i=\wt\bigwedge^5U
 \quad\hbox{and}\quad d=\wt D=s+2u,
 \]

Then the generators have weights
 \[
 \begin{array}{rclcc}
 \wt x&=& u, \\
 \wt x_{ij} &=& u+w_i+w_j, \\
 \wt x_i &=& u+w_j+w_k+w_l+w_m &=& u+s-w_i.
 \end{array}
 \]
The average of the 16 weights is $\frac12d$.

The 10 relations with their representative terms have weights
 \[
 \begin{array}{rclcl}
 N_i&=&xx_i+\cdots &\mapsto& 2u+w_j+w_k+w_l+w_m=d-w_i, \\
 N_{-i}&=&\sum x_jx_{ij} &\mapsto& d+w_i,
 \end{array}
 \]
which have average $d$. The 16 first syzygies have weights
 \[
 \begin{array}{rclcl}
 T(x) &=& \sum x_iN_{-i} &\mapsto& 2d-u, \\ 
 T(x_i) &=& xN_{-1}+\cdots &\mapsto& 2d-u-s+w_i, \\
 T(x_{ij}) &=& x_iN_j+\cdots &\mapsto& 2d-u-w_i-w_j,
 \end{array}
 \]
which have average $\frac32d$, and the 16 second syzygies
 \[
 \begin{array}{rclcl}
 S(x) &=& x^2T(x)+\cdots &\mapsto& 2d+u, \\ 
 S(x_i) &=& xx_iT(x)+\cdots &\mapsto& 2d+u+s-w_i, \\
 S(x_{ij}) &=& xx_{ij}T(x)+\cdots &\mapsto& 2d+u+w_i+w_j,
 \end{array}
 \]
which have average $\frac52s$.

To write out the Hilbert series of $\wOGr(5,10)$ with the above weights,
introduce the Laurent polynomials
 \[
 \begin{array}{rcl}
 Q_V & =& \sum_i t^{w_i}+\sum_i t^{-w_i} \\
 Q_{S^+} & =& 1+\sum_{i,j} t^{w_i+w_j}+\sum_i t^{s-w_i} \\
 Q_{S^-} & =& 1+\sum_{i,j} t^{-w_i-w_j}+\sum_i t^{-s+w_i}
 \end{array}
 \]
Then $\wOGr(5,10)$ has Hilbert series
 \[
 P(t)=1-t^{d}Q_V+t^{2d-u}Q_{S^-}-t^{2d+u}Q_{S^+}+t^{3d}Q_V-t^{4d}.
 \]
This numerology implies that the spaces of relations, first syzygies,
etc., in the resolution are the following representations of $\Spin(10)$:
 \begin{equation}
 \Oh\ot
 V\tensor D\ot
 S^- \tensor D^2\tensor L\1\ot
 S^+ \tensor D^2\tensor L\ot
 V\tensor D^3\ot
 D^4\ot 0.
 \label{eq!OGres}
 \end{equation}

Providing it is well formed, $\wOGr(5,10)$ has canonical class
$K_{\wGr}=\Oh(-4d)$. In fact the wps $\PP(S^+\tensor L)$ has
$K_{\PP}=-8d$ (the sum of weights of the coordinates), and by
(\ref{eq!OGres}), the adjunction number of $\wGr\subset\PP$ equals $4d$.

 \begin{prob}\label{prob!OGr}
 We believe that the affine O'Grassmannian $\aOGr(5,10)$ has an
equi\-variant resolution of the form (\ref{eq!OGres}), in complete analogy
with Proposition~\ref{pr!aGr}. We have written out the maps in this
sequence in explicit coordinate expressions in our treatment, with the
right $W(D_5)$ symmetry and weights. It should be possible to specify them
intrinsically in terms of Clifford multiplication.
 \end{prob}

\section{Examples}
We have searched in vain for examples of Fano 3-folds, K3 surfaces or
canonical surfaces as quasilinear sections of $\wGr(5,10)$, and we believe
that there are very few, or even none, apart from the well known straight
cases. In this section, we construct nice examples of a canonical 3-fold
and a Calabi--Yau 3-fold having isolated cyclic quotient singularities.

\begin{exa}
Let $V$ be a regular 3-fold of general type with $p_g=7$, $K^3=21$ and
$2\times\half(1,1,1)$ singularities. The plurigenus formula of Fletcher
and Reid \cite{YPG} states that
 \[
 p_n\ =\ 
 \begin{cases}
 1 \\
 p_g \\[4pt]
 \frac{n(n-1)(2n-1)}{12} K^3+(2n-1)(p_g-1)+l(n) & \hbox{for $n\ge2$,}
 \end{cases}
 \]
where $l(n)$ is a sum of the local orbifold contributions
 \[
 l(n)=
 \begin{cases}
 \frac{n}{4} \; &\hbox{if} \; n \; \hbox{is even} \\[4pt]
 \frac{n-1}{4} \;& \hbox{if} \; n \; \hbox{is odd}
 \end{cases}
 \]
 from each of the $\half(1,1,1)$ singularities. One easily calculates the
Hilbert function $H(t)=\sum p_nt^n$ from this:
 \begin{align*}
 H(t)&=1+t+\frac{t+t^2}{(1-t)^2}\cdot(p_g-1)
 +\frac{t^2+t^3}{(1-t)^4}\cdot\frac{K^3}{2}
 +\recip4\times\frac{t^2}{(1-t)(1-t^2)} \\[4pt]
 &= 1+7t+29t^2+83t^3+190t^4+370t^5+645t^6+1035t^7+1562t^8 \cdots \\[4pt]
 &= \frac{1+4t+10t^2+12t^3+10t^4+4t^5+t^6}{(1-t)^3(1-t^2)}\,.
 \end{align*}
We need seven generators in degree $1$ (since $p_g=7$) and at least two in
degree~2 to accommodate the two $\half(1,1,1)$ singularities; the simplest
possibility is that $V$ has codimension~5 in $\PP(1^7,2^2)$, with Hilbert
numerator
\[
(1-t)^7(1-t^2)^2H(t) \ =\ 1-t^2-8t^3+7t^4+8t^5-8t^7-\cdots
\]
We easily recognise this as the Hilbert numerator of the weighted
orthogonal Grassmannian $\wOGr(5,10)$ with weights $\bw=(0,0,0,0,1)$,
$u=1$ and $s=1$, therefore $d=3$. With this choice of weights, the
coordinates $x_0$, $x_{ij}$ for $1\notin \{i,j\}$ and $x_{2345}$ have
weight~1; all other coordinates have weight~2. The spinor embedding
takes $\wOGr(5,10)$ into $\PP(1^8, 2^8)$ and we construct $V$ as a
general quasilinear section
 \[
 V=\wOGr(5,10) \cap (1)\cap(2)^6.
 \]
We check that the canonical class adds up: either $V\subset\PP(1^7,2^2)$,
with adjunction number $4d=12$ gives $-7\times1-2\times2+12=1$, or
$V=(1)\cap(2)^6\subset\wOGr$ has $K_V=\Oh(-4d+1+6\times2)=\Oh(1)$.
\end{exa}

\begin{exa}
Let $V$ be a Calabi--Yau 3-fold polarised by a divisor $A$ with
$A^3=\frac65$ and $A\cdot c_2=\frac{108}{5}$, and having singular points
$P'=\recip3(1,1,1)$, $P''=\recip3(2,2,2)$, and $Q=\recip5(3,3,4)$ (we are
writing these so that $A=\Oh(1)$).

The orbifold Riemann--Roch formula of Fletcher and Reid \cite{YPG} states
that
 \[
 p_n=\frac{A^3}{6} n^3+\frac{A\cdot c_2}{12}n+c_{P'}(n)
 +c_{P''}(n)+c_{Q}(n)
 \]
where $c_\bullet(n)$ is a local contribution from the singularity that
can be calculated explicitly using the instructions in \cite{YPG}.
Following the instructions, we discover that $c_{P'}(n)+c_{P''}(n)=0$ for all
$n$, and
 \[
 c_Q(n)=
 \begin{cases}
 0 &\hbox{if $n\equiv 0$ mod 5;}\\
 0 &\hbox{if $n\equiv 1$ mod 5;}\\
 \frac{-1}{5} &\hbox{if $n\equiv 2$ mod 5;}\\
 \frac{1}{5} &\hbox{if $n\equiv 3$ mod 5;}\\
 0 &\hbox{if $n\equiv 4$ mod 5.}
 \end{cases}
 \]
{From} this it is easy to calculate the Hilbert function
 \begin{align*}
 H(t) &=
 1+\frac{A^3}6\times\frac{(1+4t+t^2)t}{(1-t)^4} +
 \frac{A\cdot c_2}{12}\times\frac{t}{(1-t)^2} +
 \recip5\times\frac{-t^2+t^3}{1-t^5} \\
 &= \frac{1-2t+3t^2-t^3-t^4+t^5+t^6-3t^7+2t^8-t^9}{(1-t)^4(1-t^5)} \\
 &=1+2t+5t^2+11t^3+20t^4+34t^5+54t^6+81t^7+117t^8+\cdots
 \end{align*}
We see that we need to multiply by $(1-t)^2(1-t^2)^2$:
 \[
 (1-t)^2(1-t^2)^2(1-t^5)H(t)=
 1+3t^3-2t^5+2t^6-3t^8-t^{11}
 \]
Then we need three generators in degree~3:
 \[
 \begin{split}
 (1-t)^2(1-t^2)^2(1-t^3)^3(1-t^5)H(t)=
 1-2t^5-4t^6+3t^8&+2t^9+\\
 &+2t^{11}+\cdots
 \end{split}
 \]
At first sight this looks like a plausible $6\times10$ codimension $4$
format; the typical example of this is a nonspecial canonical curve $C$ of
genus $6$, that is known to be a quadric section of a cone over
$\Gr(2,5)$. We might hope to find $V$ as a nonlinear section of a cone
over a weighted $\wGr(2,5)$. Indeed the polynomial in the last displayed
equation is the Hilbert numerator of
$\sC\wGr_{2,6}\cap(6)\subset\PP(1,2^3,3^6,4)$. However, this is a mirage
of a fairly typical type: although it would have the correct Hilbert
function, a quasilinear section of this variety can't have a
$\frac15(3,3,4)$ singularity. The simplest assumption is that $V$ is
codimension $5$; the easiest guess is that there is an additional
generator (and relation) in degree $4$, giving
 \begin{align*}
 & (1-t)^2(1-t^2)^2(1-t^3)^3(1-t^4)(1-t^5)H(t)=\\
 & \begin{split}
 1-t^4-2t^5-4t^6+3t^8+4t^9+4t^{10}&+2t^{11}\\
 &-2t^{13}-\cdots
 \end{split}
\end{align*}
We easily recognise this as the Hilbert numerator of the weighted
orthogonal Grassmannian $\wOGr(5,10)$ with weights $\bw=(0,0,1,1,2)$,
$u=1$ and $s=4$, therefore $d=6$, embedded in $\PP(1^2,2^4,3^4,4^4,5^2)$,
with canonical class $\Oh(-4d)=\Oh(-24)$. We can construct $V$ as a
general quasilinear section
 \[
 V=\wOGr(5,10) \cap(2)^2\cap(3)\cap(4)^3\cap(5)
 \]
(the calculation that $V$ has the correct singularities is a bit tedious
but can be done by hand).
\end{exa}

\bigskip
\noindent
Alessio Corti,\\
DPMMS, University of Cambridge,\\
Centre for Mathematical Sciences,\\
Wilberforce Road, Cambridge CB3 0WB, U.K.\\
e-mail: a.corti@dpmms.cam.ac.uk \\
web: can.dpmms.cam.ac.uk/$\!\sim$corti

\medskip
\noindent
Miles Reid,\\
Math Inst., Univ. of Warwick,\\
Coventry CV4 7AL, England\\
e-mail: miles@maths.warwick.ac.uk \\
web: www.maths.warwick.ac.uk/$\!\sim$miles

\end{document}